\documentclass[1p]{elsarticle}
\usepackage{amssymb}
\usepackage{amsfonts}

\newtheorem{theorem}{Theorem}
\newtheorem{conjecture}[theorem]{Conjecture}
\newtheorem{corollary}[theorem]{Corollary}
\newtheorem{lemma}[theorem]{Lemma}
\newtheorem{observation}[theorem]{Observation}

\newproof{pf}{Proof}

\begin{document}
\title{Decomposability of graphs into subgraphs fulfilling the 1--2--3 Conjecture}

\author[i3s]{Julien Bensmail\fnref{PEPS}}
\author[agh]{Jakub Przyby{\l}o\fnref{MNiSW}}

\fntext[PEPS]{The first author was supported by PEPS grant POCODIS.}
\fntext[MNiSW]{This work was partially supported by the Faculty of Applied Mathematics AGH UST statutory tasks within subsidy of Ministry of Science and Higher Education.}

\address[i3s]{Universit\'e C\^ote d'™Azur, CNRS, I3S, Inria, France}
\address[agh]{AGH University of Science and Technology, Faculty of Applied Mathematics, al. A. Mickiewicza 30, 30-059 Krakow, Poland}

\begin{abstract}
The well-known 1--2--3 Conjecture asserts that the edges of every graph without 
isolated edges can be weighted with $1$, $2$ and $3$ so that adjacent vertices 
receive distinct weighted degrees. This is open in general.
We prove that every $d$-regular graph, $d\geq 2$, 
can be decomposed into at most $2$ subgraphs (without isolated edges) fulfilling the 1--2--3 Conjecture if $d\notin\{10,11,12,13,15,17\}$, 
and into at most $3$ such subgraphs in the remaining cases.
Additionally, we prove that in general every graph without isolated edges can be decomposed into at most $24$ subgraphs fulfilling the 1--2--3 Conjecture, improving the previously best upper bound of $40$.
Both results are partly based on applications of the Lov\'asz Local Lemma.
\end{abstract}


\maketitle

\section{Introduction}

A graph of order at least $2$ cannot be \emph{irregular}, \emph{i.e.} its vertices cannot have 
pairwise distinct degrees. This does not concern multigraphs though. 
The least $k$ so that an irregular multigraph can be obtained from a given graph $G$ 
by replacing each edge by at most $k$ parallel edges
is called the \emph{irregularity strength} of $G$.
This graph invariant was introduced in~\cite{Chartrand}, and investigated further in numerous 
papers as a particular mean for measuring the ``level of irregularity'' of graphs, see \emph{e.g.}~\cite{Aigner,Lazebnik,Faudree2,Faudree,Frieze,KalKarPf,Lehel,MajerskiPrzybylo2,Nierhoff,Przybylo,irreg_str2}.
Potential alternative definitions of ``irregular graphs'' were also investigated by
Chartrand, Erd\H{o}s and Oellermann in~\cite{ChartrandErdosOellermann}.
In~\cite{LocalIrreg_1} the authors introduced and initiated research devoted to so-called \emph{locally irregular graphs}, \emph{i.e.} graphs in which adjacent vertices have distinct degrees.
Already earlier a local version of irregularity strength was studied in~\cite{123KLT}. 
There Karo\'nski, {\L}uczak and Thomason considered 
the least $k$ so that a locally irregular multigraph can be obtained from a given graph $G=(V,E)$ 
via, again, replacement of every edge with at most $k$ parallel edges.
This problem was in fact originally formulated in terms of weightings, where by a $k$-\emph{edge-weighting} of $G$ we mean any mapping $\omega:E\to\{1,2,\ldots,k\}$.
For such $w$ we may define the so-called \emph{weighted degree of} or simply the \emph{sum at} a given vertex $v$ as:
$$s_\omega(v):=\sum_{e\in E_v}\omega(e),$$
where $E_v$ denotes the set of edges incident with $v$ in $G$; 
we shall usually write simply $s(v)$ instead of $s_\omega(v)$ if this causes no ambiguities further on.
So in this setting, the authors of~\cite{123KLT} were interested in the least $k$ such that a $k$-edge-weighting $\omega$ of $G$ exists so that 
$s_\omega(u)\neq s_\omega(v)$ for every edge $uv\in E$ -- we say that $u$ and $v$ are \emph{sum-distinguished} then (note we must assume that $G$ contains no isolated edges to that end, \emph{i.e.} that it has no $K_2$-components).
They posed a very intriguing question, commonly known as the \emph{1--2--3 Conjecture} in the literature nowadays.
\begin{conjecture}[1--2--3 Conjecture]
For every graph $G=(V,E)$ without isolated edges there exists a weighting $\omega:E\to\{1,2,3\}$ 
sum-distinguishing all neighbours in $G$.
\end{conjecture}
They confirmed it for $3$-colourable graphs, \emph{i.e.} for graphs with $\chi(G)\leq 3$.
\begin{theorem}[\cite{123KLT}]\label{3ColoraubelFulfill123Conj}
Every $3$-colourable graph without isolated edges fulfills the 1--2--3 Conjecture.
\end{theorem}
This is also commonly known to hold in particular for complete graphs. In general the conjecture is however still widely open.
The first constant upper bound, with $30$ instead of $3$, was provided in~\cite{Louigi30}, and then improved in~\cite{Louigi} and~\cite{123with13}. The best general result thus far was delivered by Kalkowski, Karo\'nski and Pfender, who proved that it is sufficient to use weights $1,2,3,4,5$, see~\cite{KalKarPf_123}. 
This result was obtained via refinement and modification of an algorithm developed by Kalkowski~\cite{Kalkowski12} 
(
concerning a total analogue of 
the 1--2--3 Conjecture, see \emph{e.g.}~\cite{12Conjecture}).
Quite recently a complete characterization of bipartite graphs for which it is sufficient to use just weights $1$ and $2$ was also provided by Thomassen, Wu and Zhang~\cite{ThoWuZha}.
Note that graphs which require only one weight (\emph{i.e.} $1$) are precisely the locally irregular graphs.
\medskip

Another direction of research towards inducing local irregularity in a graph was developed 
by Baudon et al.~\cite{LocalIrreg_1}, this time via graph decompositions.
In this paper, by a \emph{decomposition of a graph} $G$ we mean a partition of the set of its edges into subsets inducing subgraphs of $G$ (usually with some specified features). 
We say a graph $G=(V,E)$ can be \emph{decomposed into} $k$ \emph{locally irregular subgraphs} if $E$ 
can be partitioned into $k$ sets:
$E=E_1\cup E_2\cup\ldots\cup E_k$ so that $G_i:=(V,E_i)$ is locally irregular (where we admit $E_i$ to be empty) for $i=1,2,\ldots,k$.
Equivalently, it means we may colour the edges of $G$ with at most $k$ colours so that each of these induces a locally irregular subgraph in $G$.
In~\cite{LocalIrreg_1} it was conjectured that except for some family of \emph{exceptional graphs} (each of which has maximum degree at most $3$, see~\cite{LocalIrreg_1} for details),
every connected graph can be decomposed into $3$ locally irregular subgraphs. This was then confirmed in~\cite{LocalIrreg_2} for graphs with sufficiently large minimum degree.
\begin{theorem}[\cite{LocalIrreg_2}]\label{10to10theorem}
Every graph $G$ with minimum degree $\delta(G)\geq 10^{10}$ can be decomposed into $3$ locally irregular subgraphs.
\end{theorem}
In general it was also proved by Bensmail, Merker and Thomassen~\cite{BensmailMerkerThomassen} that every connected graph which is not exceptional can be decomposed into (at most) $328$ locally irregular subgraphs,
what was then pushed down to $220$ such subgraphs by Lu\v{z}ar, Przyby{\l}o and Sot\'ak~\cite{LocalIrreg_Cubic}.
See also~\cite{LocalIrreg_1,LocalIrreg_Complexity,BensmailMerkerThomassen,LocalIrreg_Cubic} for a number of partial and related results.
\medskip

Here we develop research initiated in~\cite{8Authors}, and related  to the both concepts discussed above. 
From now on we shall write a graph \emph{fufills the 1--2--3 Conjecture} if there actually exists its neighbour sum-distinguishing $3$-edge-weighting (assuming this holds in par\-ti\-cu\-lar for an edgeless graph).
Though we are not yet able to prove the 1--2--3 Conjecture, even in the case of regular graphs,
we shall prove below that for almost every $d\geq 2$, a $d$-regular graph $G$ can be decomposed into $2$ subgraphs \emph{fulfilling the 1--2--3 Conjecture}, 
while in the remaining cases it can be decomposed into 
$3$ such subgraphs.
At the end we shall additionally prove that in general every graph without isolated edges can be decomposed into (at most) $24$ subgraphs consistent with the 1--2--3 Conjecture, while thus far it was known that $40$ such subgraphs were always sufficient, see~\cite{8Authors} (also for other related results).

\section{Basic Tools}

We first present one basic observation followed by a recollection of a few fundamental tools of 
the probabilistic method we shall use later on.  For a vertex $v$ of a given graph $G=(V,E)$, by $d_S(v)$ we mean the number of edges $uv\in E$ with $u\in S$ if $S\subseteq V$, or the number of edges $uv\in S$ in the case when $S\subseteq E$.
\begin{observation}
\label{BipartiteEvenDecomposition}
Every bipartite graph $G$ can be decomposed into two subgraphs $G_1$ and $G_2$ such that
for every vertex $v$ of $G$, 
$$d_{G_1}(v) \in \left[\frac{d_G(v)-1}{2},\frac{d_G(v)+1}{2}\right].$$
\end{observation}
\begin{pf}
If the set $U$ of the vertices of odd degree in $G$ is nonempty, then add a new vertex $u$ and join it by a single edge with every vertex in $U$; denote the obtained graph by $G'$. If $U=\emptyset$, set $G'=G$ and denote any vertex of $G'$ as $u$. As the degrees of all vertices in $G'$ are even, there exists an Eulerian tour in it. We then start at the vertex $u$ and traverse all edges of $G'$ once along this Eulerian tour colouring them alternately red and blue.
Then the red edges in $G$ induce its subgraph $G_1$ consistent with our requirements. This follows from the fact that if all degrees in the bipartite graph $G$ are even, then it has to have an even number of edges, and thus the thesis holds in particular for $u$.
\qed
\end{pf}

The following standard versions of the Lov\'asz Local Lemma can be found \emph{e.g.} in~\cite{AlonSpencer}.
\begin{theorem}[\textbf{The Local Lemma; Symmetric Version}]
\label{LLL-symmetric}
Let 
$\mathcal{A}$ be a finite family of events 
in any 
pro\-ba\-bi\-li\-ty space.
Suppose that every event $A\in \mathcal{A}$
is mutually independent of a set of all the other
events in $\mathcal{A}$ 
but at most $D$, and that 
${\rm \emph{\textbf{Pr}}}(A)\leq p$ for each $A\in \mathcal{A}$.
If
\begin{equation}\label{LLLinequality}
  ep(D+1) \leq 1,
\end{equation}
then 
${\rm \emph{\textbf{Pr}}}(\bigcap_{A\in\mathcal{A}}\overline{A})>0$.
\end{theorem}

\begin{theorem}[\textbf{The Local Lemma; General Case}]\label{LLL-general}
Let $\mathcal{A}$ be a finite family of events in any probability space and let $D=(\mathcal{A},E)$ be a directed graph such that every event $A\in \mathcal{A}$ is mutually independent of all the events $\{B: (A,B)\notin E\}$.
Suppose that there are real numbers $x_A$ ($A\in\mathcal{A}$) such that for every $A\in\mathcal{A}$, $0\leq x_A<1$ and
\begin{equation}\label{EqLLL-general}
{\rm \emph{\textbf{Pr}}}(A) \leq x_A \prod_{B\leftarrow A} (1-x_B).
\end{equation}
Then ${\rm \emph{\textbf{Pr}}}(\bigcap_{A\in\mathcal{A}}\overline{A})>0$.
\end{theorem}
Here $B\leftarrow A$ (or $A\rightarrow B$) means that there is an arc from $A$ to $B$ in $D$,
the so-called \emph{dependency digraph}.
The Chernoff Bound below can be found \emph{e.g.} in~\cite{JansonLuczakRucinski} (Th. 2.1, page 26).

\begin{theorem}[\textbf{Chernoff Bound}]\label{ChernofBoundTh}
For any $0\leq t\leq np$,
$${\mathbf Pr}({\rm BIN}(n,p)>np+t)<e^{-\frac{t^2}{3np}}~~~~{and}~~~~{\mathbf Pr}({\rm BIN}(n,p)<np-t)<e^{-\frac{t^2}{2np}}$$ 
where ${\rm BIN}(n,p)$ is the sum of $n$ independent Bernoulli variables, each equal to $1$ with probability $p$ and $0$ otherwise.
\end{theorem}

\section{Main Result for Regular Graphs}

In this section we shall prove that for almost all integers $d\geq 2$, every $d$-regular graph can be decomposed into two graphs 
fulfilling the 1--2--3 Conjecture, while in the remaining few cases -- into three such graphs. The first subsection below is devoted to small values of $d$; more generally we investigate in it 
graphs with upper-bounded chromatic number.

\subsection{Graphs with Bounded Chromatic Number}

In~\cite{8Authors} it was proved the following result (note it follows also by Corollary~\ref{ChromaticGeneralCorollary} below).
\begin{theorem}[\cite{8Authors}]\label{Chromatic9Theorem}
Every graph $G$ without isolated edges and with $\chi(G)\leq 9$ can be decomposed into $2$ graphs fulfilling the 1--2--3 Conjecture.
\end{theorem}
As complete graphs are known to fulfill the 1--2--3 Conjecture, by Brooks' Theorem we thus obtain that every $d$-regular graph with $2\leq d\leq 9$ can be decomposed into 
$2$ graphs fulfilling the 1--2--3 Conjecture.
In order to first achieve the main result of this paper for the case of regular graphs with upper-bounded degree, we
shall generalize Theorem~\ref{Chromatic9Theorem}. 
For this aim we 
use the following Lemma~\ref{2-colour_lemma}, which is very similar 
to a one proved in~\cite{8Authors}. 
We however present an alternative brief proof of this slightly modified version
for the sake of completeness of the exposition of our reasoning. 

\begin{lemma}
\label{2-colour_lemma}
If the edges of a graph $G$ without isolated edges 
can be $2$-coloured with red and blue so that the induced red subgraph $R$ and blue subgraph $B$ satisfy $\chi(R)\leq r\geq 3$ and $\chi(B)\leq b\geq 3$, then we can also do it in such a way that neither $B$ nor $R$ contains an isolated edge.
\end{lemma}

\begin{pf}
Start from a $2$-colouring of the edges of $G$ with $\chi(R)\leq r\geq 3$ and $\chi(B)\leq b\geq 3$ which minimizes the number of monochromatic $K_2$-components
(note that each isolated triangle 
of $G$ shall be monochromatic then). We shall show that if the number of these is still positive, then we may ``get rid'' of any given such component, without creating a new one (thus getting a contradiction, and hence proving the thesis):

Assume $uv$ forms such a monochromatic, say blue $K_2$-component. Observe that:

(1') $u$ and $v$ must belong to the same component of $R$, as otherwise one may recolour $uv$ red;

(2') if $e$ is a red edge adjacent with $uv$, say $e=uw$, then the size $k$ of the (red) component of $R-e$ including $w$ equals $1$ (hence, all red paths originating at $u$ or $v$ must be of length $2$, \emph{i.e.} there are in particular no isolated red edges adjacent with $uv$), as otherwise, if there was any such edge $e$ with $k=0$, we could recolour $uw$ blue (we would not create any red $K_2$-component then, as we would not change colours on some existing red path joining $u$ and $v$), while in the remaining cases (\emph{i.e.} when $k\geq 2$ for each red edge $e$ adjacent with $uv$) we could also recolour $uw$ blue. 

By (1') and (2') there must be a red path $uwv$ in  $G$, and hence, by (2'), $d(u)=2=d(v)$. Therefore, as all 
isolated triangles of 
$G$ are monochromatic, 
$d(w)\geq 3$, and thus by (2') all edges incident with $w$ except for 
$uw$ and $vw$ are blue. We may however recolour $uv$ red and $uw$ blue then.
\qed
\end{pf}

\begin{lemma}\label{3-colourable_decomposition}
For each positive integer $k$, every graph $G=(V,E)$ without isolated edges and with $\chi(G)\leq 3^k$ can be decomposed into $k$ (some possibly empty) subgraphs $G_1,G_2,$ $\ldots,$ $G_k$ such that $\chi(G_i)\leq 3$ and $G_i$ contains no isolated edges for $i=1,2,\ldots,k$.
\end{lemma}

\begin{pf}
We prove the lemma by induction with respect to $k$.
For $k=1$ it trivially holds, so let us assume that $k\geq 2$.
Partition $V$ into (possibly empty) independent sets $V_1,V_2,\ldots,V_{3^k}$.
Colour red every edge $uv\in E$ such that $u\in V_i$ and $v\in V_j$ with $i\equiv j~({\rm mod}~3)$, and colour blue the remaining edges of $G$.
Let $R$ and $B$ be the red and, resp., blue subgraphs of $G$. Note that $\chi(B)\leq 3$, as $V_{1+r}\cup V_{4+r}\cup V_{7+r}\cup\ldots\cup V_{3^k-2+r}$ forms an independent set in $B$ for $r=0,1,2$.
On the other hand, $\chi(R)\leq 3^{k-1}$, as there are no edges in $R$ between any two of the three red subgraphs 
$G[V_{1+r}\cup V_{4+r}\cup V_{7+r}\cup\ldots\cup V_{3^k-2+r}]$ with $r=0,1,2$.
Note that by Lemma~\ref{2-colour_lemma}, we may assume that neither $B$ nor $R$ contains an isolated edge.
By the induction hypothesis, $R$ can be then decomposed into subgraphs $G_1,G_2,\ldots,G_{k-1}$ with $\chi(G_i)\leq 3$ none of which contains an isolated edge. Setting $G_k:=B$ thus yields the thesis.
\qed
\end{pf}

By Lemma~\ref{3-colourable_decomposition} and Theorem~\ref{3ColoraubelFulfill123Conj} we obtain the following corollaries.

\begin{corollary}\label{ChromaticGeneralCorollary}
Every graph $G$ without isolated edges can be decomposed into $\left\lceil\log_3\chi(G)\right\rceil$ graphs fulfilling the 1--2--3 Conjecture.
\end{corollary}

\begin{corollary}\label{ChromaticRegularCorollary}
Every $d$-regular graph $G$ with $10\leq d\leq 17$ can be decomposed into $3$ subgraphs fulfilling the 1--2--3 Conjecture.
\end{corollary}

To prove divisibility into 2 such subgraphs 
of $d$-regular graphs with larger $d$, we first prove in the next Subsection~\ref{SubsectRandVertPart} that they admit a special vertex partition. In Subsection~\ref{Subsection123Family}, we then discuss
a peculiar sufficient condition for a graph to fulfill the 1--2--3 Conjecture. 

\subsection{Random Vertex Partition}\label{SubsectRandVertPart}

\begin{lemma}\label{ProbabilisticPartition23}
The vertices of every $d$-regular graph $G$ with $d\geq 14$, $d\neq 15,17$, can be partitioned into sets $V_0$ and $V_1$ such that if $d\equiv r~{\rm mod}~2$ for some $r\in\{0,1\}$, then:
\begin{itemize}
\item[(i)] $\forall v\in V_0: d_{V_0}(v)\geq 2+r$;
\item[(ii)] $\forall v\in V_0: d_{V_1}(v)\geq 2$;
\item[(iii)] $\forall v\in V_1: d_{V_1}(v)\geq 2+r$;
\item[(iv)] $\forall v\in V_1: d_{V_0}(v)\geq 2$.
\end{itemize}
\end{lemma}

\begin{pf}
Assume $G=(V,E)$ is a $d$-regular graph with $d\geq 14$, $d\neq 15,17$.
To every vertex we randomly and independently assign $0$ or $1$ -- each with probability $1/2$,
and denote the 2-colouring obtained by $c$. Set $V_0=c^{-1}(0), V_1=c^{-1}(1)$.

Assume first that $d$ is even (\emph{i.e.}, $r=0$). For any given vertex $v\in V$, denote by:
\begin{itemize}
\item $A_1(v)$ -- the event that $v\in V_0$ and $d_{V_0}(v)\leq 1$;
\item $A_2(v)$ -- the event that $v\in V_0$ and $d_{V_1}(v)\leq 1$;
\item $A_3(v)$ -- the event that $v\in V_1$ and $d_{V_1}(v)\leq 1$;
\item $A_4(v)$ -- the event that $v\in V_1$ and $d_{V_0}(v)\leq 1$.
\end{itemize}
Note that if we prove that none of these events holds for some (random) colouring $c$, then the thesis shall be fulfilled. As drawings for all vertices are independent, for every $v\in V$ we have: 
\begin{eqnarray}\label{EvenEquation1}
{\mathbf Pr}(A_1(v)) &=& {\mathbf Pr}(v\in V_0) \cdot {\mathbf Pr}(d_{V_0}(v)\leq 1) \nonumber\\
&=& {\mathbf Pr}(v\in V_0) \cdot \left[{\mathbf Pr}(d_{V_0}(v) = 0) + {\mathbf Pr}(d_{V_0}(v) = 1)\right] \nonumber\\
&=& \frac{1}{2}\cdot \left[\left(\frac{1}{2}\right)^d+d\left(\frac{1}{2}\right)^d\right] = (1+d)\left(\frac{1}{2}\right)^{d+1},
\end{eqnarray}
and analogously, for every $i=2,3,4$,
\begin{equation}\label{EvenEquations234}
{\mathbf Pr}(A_i(v)) = (1+d)\left(\frac{1}{2}\right)^{d+1}.
\end{equation}
Note now that every event $A_i(v)$ is mutually independent of all other events $A_j(u)$ with $u$ at distance at least $3$ from $v$, \emph{i.e.} of all but at most $4d^2+3$ other events. Therefore, by (\ref{LLLinequality}), (\ref{EvenEquation1}) and~(\ref{EvenEquations234}), in order to apply the Lov\'asz Local Lemma 
it is 
sufficient to show that:
\begin{equation}\label{LLLappl1}
  e(1+d)\left(\frac{1}{2}\right)^{d+1}(4d^2+4)<1
\end{equation}
for every even integer $d\geq 14$. For $d=14$ the left-hand side of inequality~(\ref{LLLappl1}) takes value (approximately) $0.9805...<1$, while for $d\geq 16$ inequality~(\ref{LLLappl1}) is implied by the following one:
$$2e(1+d)^3<2^d,$$
equivalent to:
$$\sqrt[3]{2e}(1+d)<2^{\frac{d}{3}},$$
which holds as for $f(d):=2^{d/3}-\sqrt[3]{2e}(1+d)$, we have $f'(d)=2^{d/3} 1/3 \ln2 -\sqrt[3]{2e}>0$ for $d\geq 16$ and $f(16)\approx 10.4253>0$.

By (\ref{EvenEquation1}), (\ref{EvenEquations234}), (\ref{LLLappl1}) and Theorem~\ref{LLL-symmetric} 
we thus conclude that
$${\mathbf Pr}\left(\bigcap_{v\in V} (\overline{A_1(v)}\cap \overline{A_2(v)}\cap \overline{A_3(v)}\cap \overline{A_4(v)})\right)>0.$$
The thesis follows.\medskip

Assume now that $d$ is odd (\emph{i.e.}, $r=1$) and $d\geq 19$.
In order to optimize our approach we shall now have to
aggregate the events concerning our requirements \emph{(i)--(iv)}.
Thus for a vertex $v\in V$, denote the following (aggregated) event:
\begin{itemize}
\item $B(v)$: $(v\in V_0 \wedge d_{V_0}(v)\geq 3 \wedge d_{V_1}(v)\geq 2) \vee (v\in V_1 \wedge d_{V_1}(v)\geq 3 \wedge d_{V_0}(v)\geq 2)$.
\end{itemize}
In order to prove the thesis it is then enough to apply the Local Lemma to show that the probability that $B(v)$ holds for every $v\in V$ is positive.

Note that for every $v\in V$, 
\begin{eqnarray}\label{OddEquation0}
{\mathbf Pr}\left(\overline{B(v)}\right) &=& {\mathbf Pr}\left((v\in V_1 \vee d_{V_0}(v)\leq 2 \vee d_{V_1}(v)\leq 1) \wedge (v\in V_0 \vee d_{V_1}(v)\leq 2 \vee d_{V_0}(v)\leq 1)\right) \nonumber\\
 &=& \frac{1}{2}{\mathbf Pr}\left(\overline{B(v)}|v\in V_1\right) + \frac{1}{2}{\mathbf Pr}\left(\overline{B(v)}|v\in V_0\right)  \nonumber\\
 &=& {\mathbf Pr}\left(\overline{B(v)}|v\in V_1\right) \nonumber\\
 &=& {\mathbf Pr}\left(d_{V_1}(v)\leq 2 \vee d_{V_0}(v)\leq 1\right)  \nonumber\\
 &=& {\mathbf Pr}\left(d_{V_1}(v) = 2\right) + {\mathbf Pr}\left(d_{V_1}(v) = 1\right) + {\mathbf Pr}\left(d_{V_1}(v) = 0\right) \nonumber\\
 &&+ {\mathbf Pr}\left(d_{V_0}(v) = 1\right) + {\mathbf Pr}\left(d_{V_0}(v) = 0\right) \nonumber\\
 &=& {d\choose 2}\left(\frac{1}{2}\right)^d + d\left(\frac{1}{2}\right)^d + \left(\frac{1}{2}\right)^d + d\left(\frac{1}{2}\right)^d + \left(\frac{1}{2}\right)^d \nonumber\\
 &=& \left(\frac{d(d-1)}{2}+2d+2\right)\left(\frac{1}{2}\right)^d \nonumber\\
 &<& (d+2)^2 2^{-(d+1)}.
\end{eqnarray}
Again an event $\overline{B(v)}$ is mutually independent of all other events $\overline{B(u)}$ with $u$ at distance at least $3$ from $v$, \emph{i.e.} of all but at most $d^2< (d+2)^2-1$ other events.
Moreover, the following inequality:
\begin{equation}\label{LLLappl2}
  e (d+2)^2 2^{-(d+1)}
  (d+2)^2<1
\end{equation}
is fulfilled for every $d\geq 19$, as this is equivalent to the 
fact that $g(d):=2^{(d+1)/4}-\sqrt[4]{e}(d+2)> 0$ (for $d\geq 19$). 
This in turn holds, since $g(19)\approx 5.0355>0$ and $g'(d)= 2^{(d-7)/4}\ln 2 -\sqrt[4]{e}>0$ for $d\geq 19$.

By~(\ref{OddEquation0}), (\ref{LLLappl2})  and Theorem~\ref{LLL-symmetric}
we thus obtain that 
$${\mathbf Pr}\left(\bigcap_{v\in V} B(v)\right)>0.$$
\qed
\end{pf}

\subsection{Family of Graphs Fulfilling the 
1--2--3 Conjecture}\label{Subsection123Family}

For a subset $S$ of the set of vertices of a given graph $G$, 
$G[S]$ 
shall denote the subgraph induced by $S$ in $G$, 
while by $G_1\cup G_2$ we shall mean the sum of two graphs $G_1=(V_1,E_1)$, $G_2=(V_2,E_2)$ understood as the pair $(V_1\cup V_2,E_1\cup E_2)$.
An \emph{independent set} in a graph $G=(V,E)$ is a subset $I$ of $V$ such that no edge of $G$ has both ends in $I$.
We call it \emph{maximal} (or an \emph{independent dominating set}) if every vertex in $V\smallsetminus I$ has a neighbour in $I$. In order to prove the existence of a specific family of graphs fulfilling the 1--2--3 Conjecture we shall apply in the following lemma a certain refinement of Kalkowki's algorithm from~\cite{Kalkowski12}, exploiting 
for this aim
the concept of maximal independent sets; see~\cite{Julien5regular123} for a corresponding application of a mixture of these two ingredients.

\begin{lemma}\label{123ConjGraphFamily}
If a graph $G=(V,E)$ contains a maximal independent set $I$ such that there exists a constant $\alpha\geq 1$ so that for $R:=V\smallsetminus I$,
\begin{itemize}
\item[$(1^\circ)$] $d(v)\leq \alpha$ for every $v\in I$ and 
\item[$(2^\circ)$] $d(v)\geq \alpha+\frac{d_R(v)+1}{2}$ for every $v\in R$, 
\end{itemize}
then $G$ fulfills the 1--2--3 Conjecture.
\end{lemma}

\begin{pf} By $(2^\circ)$, $d(v)\geq 2$ for every $v\in R$, and thus there are no isolated edges 
in $G$.
Note also that since $I$ is a maximal independent set in $G$, then:
\begin{itemize}
\item[$(3^\circ)$] $d_I(v)\geq 1$ for every $v\in R$,
\end{itemize}
so for every vertex $v\in R$ we may fix an edge $e_v$ joining $v$ with some vertex in $I$.

We shall construct a $3$-edge-weighting of $G$ sum-distinguishing its neighbours.
Initially we label all edges in $G$ by $2$. These shall be modified gradually, and by $\omega(e)$ we shall always understand 
the current weight of an edge $e$ in a given moment of our on-going relabelling algorithm specified below, and similarly, by $s(v)$ we shall understand 
the current sum at a vertex $v$ in $G$.

Let $G_R=G[R]$ be the graph induced by $R$ in $G$. We analyse every of its components one by one, in any fixed order, and modify the labels of some of the edges incident (in $G$) with at least one vertex of this component. Suppose $H$ is the next component to be analysed within the on-going algorithm,
and arbitrarily order its vertices linearly into a sequence $v_1,v_2,\ldots,v_n$. We shall analyse one vertex in the sequence after another starting from $v_1$, for which we perform no changes. Suppose thus we are about to analyse a vertex $v_j$ with $j\geq 2$ (if $H$ has more than one vertex) and denote by $N_H^-(v_j)$ the set of neighbours of $v_j$ in $H$ which precede it in the fixed linear ordering -- we call these the \emph{backward neighbours} of $v_j$.
Similarly we define the set $E_H^-(v_j)$ of the \emph{backward edges} of $v_j$, \emph{i.e.} these joining $v_j$ with its backward neighbours in $H$.
We shall now modify (if necessary) weights of some edges, in order to obtain a sum at $v_j$ (in $G$) which is distinct from the sums of all its backward neighbours -- this sum of $v$ shall then not change in the further part of the algorithm. To obtain our goal we shall be allowed to perform changes only on the edges incident with its backward neighbours, namely for every backward edge $v_kv_j\in E_H^-(v_j)$ of $v_j$ (\emph{i.e.} with $k<j$) we shall be allowed to modify the labels of $v_kv_j$ and $e_{v_k}$ so that the sum at $v_k$ does not change; more specifically, if prior to this step $j$ we had $\omega(e_{v_k})=2$ (and $\omega(v_kv_j)=2$), then we may increase the label of $v_kv_j$ by $1$ and decrease the label of $e_{v_k}$ by $1$ (or perform no changes on these two edges), while if priory we had $\omega(e_{v_k})=1$, we may decrease the label of $v_kv_j$ by $1$ and increase the label of $e_{v_k}$ by $1$ (observe that then $\omega(e_{v_k})\in\{1,2\}$ and $\omega(v_kv_j)\in\{1,2,3\}$). Note that such admitted operations allow us to change the label of every backward edge of $v_j$ by exactly $1$ (or do nothing with this label). Hence we have available at least $|E_H^-(v_j)|+1$ distinct sums at $v_j$ via these operations. We choose one of these sums which is distinct from the current sums of all backward neighbours of $v_j$ and denote it by $s^*$ (it exists, as $|E_H^-(v_j)|+1=|N_H^-(v_j)|+1$), and we perform (some of the) admissible changes described above so that $s(v_j)=s^*$ afterwards. As these changes do not influence sums (in $G$) of the other vertices of $H$, 
$v_j$ is now sum-distinguished from all its backward neighbours, and  
\begin{eqnarray}
s(v_j) &\geq& 2\cdot d_I(v_j)+1\cdot d_R(v_j) = 2(d(v_j)-d_R(v_j))+d_R(v_j) \nonumber\\
       &=& 2\left(d(v_j)-\frac{d_R(v_j)}{2}\right)\geq 2\alpha+1 \label{DegreIneqR}
\end{eqnarray}
by $(2^\circ)$. This shall not change as we guarantee that the sums of $v_1,v_2,\ldots,v_j$ shall not be modified in the further part of the construction.
After step $n$, all neighbours in $H$ are thus sum-distinguished (in $G$), and we continue in the same manner with a consecutive component of $G_R$, if any is still left. Note that performing such changes concerning one component of $G_R$ does not influence the sums in the other components, hence at the end of our construction all neighbours in $G_R$ are sum-distinguished (in $G$).
On the other hand, as due to our algorithm every edge incident with a vertex in $I$ has final weight $1$ or $2$, by $(1^\circ)$ we obtain that for $v\in I$, 
\begin{equation}\label{DegreIneqI}
s(v)\leq 2d(v)\leq 2\alpha.
\end{equation}
Hence, by~(\ref{DegreIneqR}) and~(\ref{DegreIneqI}), every vertex in $R$ is also sum-distinguished from each of its neighbours in $I$.
As $I$ is an independent set, we thus obtain a desired $3$-edge-weighting of $G$.
\qed
\end{pf}

\subsection{Main Result for Almost All Degrees}

\begin{theorem}\label{MainRegularDecomposition}
Every $d$-regular graph $G$ with $d\geq 14$, $d\neq 15,17$, can be decomposed into two graphs fulfilling the 1--2--3 Conjecture.
\end{theorem}

\begin{pf}
Let $G=(V,E)$ be a $d$-regular graph with $d\geq 14$, $d\neq 15,17$, and
$d\equiv r~{\rm mod}~2$ for some $r\in\{0,1\}$.
Let then $V=V_0\cup V_1$ be a vertex partition consistent with the thesis of Lemma~\ref{ProbabilisticPartition23}.
Denote $G_0:=G[V_0]$, $G_1:=G[V_1]$.

Let $H$ be the bipartite graph induced by the edges between $V_0$ and $V_1$. Then $\delta(H)\geq 2$ by \emph{(ii)} and \emph{(iv)} from Lemma~\ref{ProbabilisticPartition23}.
By Observation~\ref{BipartiteEvenDecomposition} we next decompose $H$ into two subgraphs $H_0$ and $H_1$ such that
\begin{equation}\label{H0degrees}
d_{H_0}(v) \in \left[\frac{d_H(v)-1}{2},\frac{d_H(v)+1}{2}\right],
\end{equation}
and thus also
\begin{equation}\label{H1degrees}
d_{H_1}(v) \in \left[\frac{d_H(v)-1}{2},\frac{d_H(v)+1}{2}\right]
\end{equation}
for every vertex $v\in V$.

Let $G'_0=G_0\cup H_0$, $G'_1=G_1\cup H_1$. Obviously $G'_0$ and $G'_1$ constitute a decomposition of $G$.
In order to finish the proof it is thus sufficient to prove that they are both consistent with the assumptions of Lemma~\ref{123ConjGraphFamily}.
We show this to hold for $G'_0$, as the reasoning for $G'_1$ is precisely symmetrical.

For this aim note first that by the definition of $H_0$, the set $I:=V_1$ is an independent set in $G'_0$,
and by \emph{(ii)} from Lemma~\ref{ProbabilisticPartition23} and~(\ref{H0degrees}) above, it is also maximal. We shall now show that $(1^\circ)$ and $(2^\circ)$ from Lemma~\ref{123ConjGraphFamily} hold for $G'_0$ with $R=V_0$ and $\alpha:=\frac{d-2-r}{2}$.

By \emph{(iii)} from Lemma~\ref{ProbabilisticPartition23}, for every vertex $v\in I=V_1$,
$$d_H(v)\leq d-2-r,$$
and hence, by~(\ref{H0degrees}): 
$$d_{G'_0}(v)=d_{H_0}(v) \leq \left\lceil\frac{d_H(v)}{2}\right\rceil \leq \left\lceil \frac{d-2-r}{2}\right\rceil = \frac{d-2-r}{2}$$
according to the definition of $r$. Consequently, $(1^\circ)$ holds.

On the other hand, by~(\ref{H0degrees}), for every $v\in R = V_0$:
\begin{eqnarray}
d_{G'_0}(v) &=& d_{G_0}(v) + d_{H_0}(v) \geq d_{V_0}(v)+\left\lfloor\frac{d-d_{V_0}(v)}{2}\right\rfloor
\geq d_{V_0}(v)+\frac{d-d_{V_0}(v)-1}{2} \nonumber\\
&=& \frac{d+d_{V_0}(v)-1}{2} \geq \frac{d+d_{V_0}(v)-1-r}{2} = \frac{d-2-r}{2} + \frac{d_{V_0}(v)+1}{2} \nonumber\\
&=& \alpha + \frac{d_{R}(v)+1}{2}, \nonumber
\end{eqnarray}
and thus $(2^\circ)$ holds.
\qed
\end{pf}

\section{General Upper Bound for All Graphs}

We conclude by showing that every graph without isolated edges can be decomposed into a certain number $K$ of  graphs fulfilling the 1--2--3 Conjecture,
where $K\leq 24$. 
We thereby improve the previously best upper bound $K\leq 40$ from~\cite{8Authors}.
We start from proving a lemma on the existence of a subset of edges with certain properties in graphs with sufficiently large minimum degree.

\begin{lemma}\label{Edge_subset_S_lemma}
If $G=(V,E)$ is a graph with minimum degree $\delta \geq 10^{10}+10^8$, 
then there is a subset $S\subseteq E$ such that
$1\leq d_S(v)\leq d(v)-10^{10}$  for every vertex $v\in V$.
\end{lemma}

\begin{pf}
Let $\Delta=\Delta(G)$.
For every vertex $v\in V$ choose arbitrarily a subset $F_v\subseteq E_v$ of cardinality $\delta$.
Now randomly and independently for every vertex $v\in V$ 
choose one edge in $F_v$ -- each with equal probability (\emph{i.e.} $\delta^{-1}$) -- and denote it by $e_v$.
For every $v\in V$ denote the event:
\begin{itemize}
\item $A(v)$: $|\{u\in N_{G}(v): e_u\neq uv\}| < 10^{10}+1$.
\end{itemize}
Suppose $v$ is a vertex of degree $d$; note that by the Chernoff Bound:
\begin{eqnarray}\label{Probabil_subset_v}
{\mathbf Pr}\left(A(v)\right) &\leq& {\mathbf Pr}\left({\rm BIN}\left(d,\frac{\delta-1}{\delta}\right) < 10^{10}+1\right) 
< e^{-\frac{\left(\frac{\delta-1}{\delta}d-10^{10}-1\right)^2}{2\frac{\delta-1}{\delta}d}}\nonumber\\
&\leq& e^{-\frac{\left(\frac{\delta-1}{\delta}d-\frac{\delta-1}{\delta}d\left(1-10^{-3}\right)\right)^2}{2\frac{\delta-1}{\delta}d}} <  e^{-2\cdot 10^{-7}d}, 
\end{eqnarray}
and set $x_v:=e^{-10^{-7}d}$ for such $v$. In order to apply the general version of the Local Lemma we define a dependency digraph $D$ by joining $A(v)$ with an arc to every $A(u)$ for which there exists $w\in V$ such that  $uw,vw\in F_w$; note that there are at most $d\delta$ such events $A(u)$ for $v$. Then, since $e^{-x}<1-x+0,5x^2$ for $x>0$ and $f(x):=e^{-10^{-8}x}x$ is decreasing for $x\geq 10^8$,
we have:
\begin{eqnarray}\label{Probabil_subset_v_LLL}
x_v \prod_{A(u)\leftarrow A(v)} (1-x_u) &\geq& e^{-10^{-7}d} \left(1-e^{-10^{-7}\delta}\right)^{\delta d} \nonumber\\
&>& e^{-10^{-7}d} \left(1-e^{-10^{-8}\delta}+0.5e^{-2\cdot 10^{-8}\delta}\right)^{\delta d} \nonumber\\
&>& e^{-10^{-7}d} \left(e^{-e^{-10^{-8}\delta}}\right)^{\delta d} \nonumber\\
&>& e^{-10^{-7}d} \left(e^{-e^{-10^{-8}10^{10}}10^{10}}\right)^d \nonumber\\
&>& e^{-2\cdot 10^{-7}d}.
\end{eqnarray}
By (\ref{Probabil_subset_v}), (\ref{Probabil_subset_v_LLL}) and Theorem~\ref{LLL-general} we thus 
conclude that there is a choice of edges $e_u$, $u\in V$, so 
that for every $v\in V$, 
$$|\{u\in N_{G}(v): e_u\neq uv\}| \geq 10^{10}+1.$$
It is then sufficient to set $S=\{e_u:u\in V\}$ to obtain $1\leq d_S(v)\leq d(v)-10^{10}$ for each $v\in V$, as desired.
\qed
\end{pf}

We are now ready to prove a lemma resembling one of the observations (\emph{i.e.} Lemma~4.5) from~\cite{BensmailMerkerThomassen} (used there as an ingredient in research concerning graph decompositions into a given finite number of locally irregular subgraphs). Lemma~\ref{Edge_subset_S_lemma} above shall enable us to optimize 
the thesis of the aforementioned 
Lemma~\ref{degenarated_min_deg_decomp} below.

\begin{lemma}\label{degenarated_min_deg_decomp}
Every graph $G=(V,E)$ without isolated edges can be decomposed into
two graphs $H$ and $F$ such that: 
$H$ is either empty or has minimum degree $\delta(H)\geq 10^{10}$, and 
$F$ contains no isolated edges and has degeneracy less than $10^{10}+10^8$.
\end{lemma}

\begin{pf}
We shall first gradually remove some vertices from a given graph $G$. As long as there is still some vertex $v$ of degree less than $10^{10}+10^8$ in what is left of it, we remove $v$ from our contemporary graph. At the end of this process, we denote the leftover of $G$ by $H'$ and let $F'$ be the subgraph of $G$ induced by all its edges with at least one end outside $V(H')$. Note that $F'$ has degeneracy less than $10^{10}+10^8$, while $\delta(H')\geq 10^{10}+10^8$ or $H'$ is empty then.
If $H'$ is empty, the thesis holds. Otherwise, by Lemma~\ref{Edge_subset_S_lemma}, there exists
$S\subseteq E(H')$ such that for every vertex $v\in V(H')$,
\begin{equation}\label{Inequalities_from_ESSL}
1\leq d_S(v)\leq d_{H'}(v)-10^{10}.
\end{equation}
Note that for every isolated edge $uv$ of $F'$, one of its ends must belong to $V(H')$ -- we then arbitrarily choose one edge from $S$ incident with this end and add it to $F'$ provided that no other edge adjacent to $uv$ was earlier added to $F'$. After repeating this procedure for every such isolated edge we obtain a graph $F$ of $F'$; note that
the degeneracy of $F$ is still less than $10^{10}+10^8$ (as we may place 
the ends of the isolated edges of $F'$
together with the vertices in $V(F)\smallsetminus V(F')$ at the end of the ordering witnessing the degeneracy of $F$, since these vertices induce a forest in $F$). At the same time, by~(\ref{Inequalities_from_ESSL}), the remaining subgraph of $G$, denoted by $H$ (formed from $H'$ by removing edges from $S$ transferred to $F'$), fulfills: $\delta(H)\geq 10^{10}$.
\qed
\end{pf}

\begin{theorem}
Every graph $G$ without isolated edges can be decomposed into $24$ subgraphs fulfilling the 1--2--3 Conjecture.
\end{theorem}

\begin{pf}
By Lemma~\ref{degenarated_min_deg_decomp}, $G$ can be decomposed into a graph $H$ which is either empty or has minimum degree $\delta(H)\geq 10^{10}$ and
a graph $F$ of degeneracy less than $10^{10}+10^8$ which contains no isolated edges. By Theorem~\ref{10to10theorem}, $H$ can be further decomposed into $3$ 
locally irregular subgraphs (which obviously fulfill the 1--2--3 Conjecture).
On the other hand, as $\chi(F)\leq 10^{10}+10^8 < 3^{21}$, by Lemma~\ref{3-colourable_decomposition}, $F$ can be decomposed into 
$21$ graphs which are $3$-colourable
and contain no isolated edges,
and thus fulfill the 1--2--3 Conjecture by Theorem~\ref{3ColoraubelFulfill123Conj}.
\qed
\end{pf}

\section{Concluding Remarks}

Note that by Theorems~\ref{Chromatic9Theorem} and~\ref{MainRegularDecomposition} we know that any $d$-regular graph without isolated edges can be decomposed into $2$ subgraphs fulfilling the 1--2--3 Conjecture if only $d\notin \{10,11,12,13,15,17\}$. The remaining cases apparently need a separate special treatment, but either way, by Corollary~\ref{ChromaticRegularCorollary} every $d$-regular graph, $d\geq 2$, can be decomposed into (at most) $3$ subgraphs complying with the 1--2--3 Conjecture.
An even more challenging task was proposed in~\cite{8Authors}, whose authors suspect that something stronger should hold. Namely, they conjectured that in fact every graph $G$ without isolated edges and isolated triangles can be decomposed into $2$ subgraphs fulfilling the 1--2--3 Conjecture with only weights $1$ and $2$
(\emph{i.e.}, admitting neighbour sum-distinguishing 2-edge-weightings).
This interesting problem is independent of the 1--2--3 Conjecture itself, and is also partly related to the research from~\cite{BensmailMerkerThomassen} -- it is in particular known that this conjecture holds for bipartite graphs and subcubic graphs.
See~\cite{8Authors} for details and further observations 
concerning this new concept, and many other related problems and results.

\end{document}